\documentclass[12pt]{amsart}
\usepackage{amsfonts,amssymb}
\RequirePackage{ifpdf}
\usepackage{amsmath}
\usepackage{color}
\usepackage{pstcol,pst-node}
\usepackage{graphics}
\usepackage{epic}
\usepackage{curves}

\def\appendix#1{
\addtocounter{section}{1} \setcounter{equation}{0}
\renewcommand{\thesection}{\Alph{section}}
\section*{Appendix \thesection\protect\indent\quad
#1}
}
\renewcommand{\theequation}{\thesection.\arabic{equation}}

\catcode`\@=11
\def\marginnote#1{}

\newcount\hour
\newcount\minute
\newtoks\amorpm
\hour=\time\divide\hour by60 \minute=\time{\multiply\hour by60
\global\advance\minute by-\hour}
\edef\standardtime{{\ifnum\hour<12 \global\amorpm={am}%
        \else\global\amorpm={pm}\advance\hour by-12 \fi
        \ifnum\hour=0 \hour=12 \fi
        \number\hour:\ifnum\minute<10 0\fi\number\minute\the\amorpm}}
\edef\militarytime{\number\hour:\ifnum\minute<100\fi\number\minute}

\def\draftlabel#1{{\@bsphack\if@filesw {\let\thepage\relax
      \xdef\@gtempa{\write\@auxout{\string
          \newlabel{#1}{{\@currentlabel}{\thepage}}}}}\@gtempa \if@nobreak
    \ifvmode\nobreak\fi\fi\fi\@esphack} \gdef\@eqnlabel{#1}}
    \def\@eqnlabel{}
\def\@vacuum{}
\def\draftmarginnote#1{\marginpar{\raggedright\scriptsize\tt#1}}

\def\draft{
  \oddsidemargin -.5truein
  \def\@oddfoot{\footnotesize \sl preliminary draft \hfil
    \rm\thepage\hfil\sl\today\quad\militarytime}
  \let\@evenfoot\@oddfoot \overfullrule 3pt
    \let\label=\draftlabel
    \let\marginnote=\draftmarginnote
  \def\@eqnnum{(\theequation)\rlap{\kern\marginparsep\tt\@eqnlabel}%
    \global\let\@eqnlabel\@vacuum}

  }

\voffset= - 0.5in \hoffset= - 1.0in

\newcommand{\tr}{\,{\rm Tr}\,}

\def\be{\begin{equation}}
\def\ee{\end{equation}}
\def\bea{\begin{eqnarray}}
\def\eea{\end{eqnarray}}
\def\<{\langle}
\def\>{\rangle}

\def\tr{\mathop{\rm{tr}}}

\def\ocomma{{\phantom{\Bigm|}^{\phantom {X}}_{\raise-1.5pt\hbox{,}}\!\!\!\!\!\!\otimes}}

\newtheorem{theorem}{Theorem}[section]
\newtheorem{lm}[theorem]{Lemma}

\theoremstyle{definition}

\newtheorem{remark}[theorem]{Remark}

\allowdisplaybreaks

\begin{document}

\title[The matrix model for dessins d'enfants]
{The matrix model for dessins d'enfants}
\author{Jan Ambj{\o}rn$^\dagger$}\thanks{$^\dagger$Niels Bohr Institute,
Copenhagen University, Denmark, and IMAPP, Radboud University, Nijmengen, The Netherlands.
Email: ambjorn@nbi.dk.}
\author{Leonid Chekhov$^{\ast}$}\thanks{$^{\ast}$Steklov Mathematical
Institute and  Laboratoire Poncelet, Moscow, Russia; Niels Bohr Institute, Copenhagen University, Denmark. Email:
chekhov@mi.ras.ru.}

\begin{abstract}
We present the matrix models that are the generating functions for branched
covers of the complex projective line ramified over $0$, $1$, and $\infty$
(Grotendieck's dessins d'enfants) of fixed genus, degree, and the ramification
profile at infinity. For general ramifications at other points, the model
is the two-logarithm matrix model with the external field studied previously by
one of the authors (L.Ch.) and K.Palamarchuk. It lies in the class of the
generalised Kontsevich models (GKM) thus being the Kadomtsev--Petviashvili (KP) hierarchy
$\tau$-function and, upon the shift of times, this model is equivalent to
a Hermitian one-matrix model with a general potential whose coefficients are
related to the KP times by a Miwa-type transformation.
The original model therefore enjoys a topological recursion and can be
solved in terms of shifted moments of the standard Hermitian one-matrix model at all genera
of the topological expansion. We also derive the matrix model for clean Belyi morphisms, which
turns out to be the Kontsevich--Penner model introduced by the authors and Yu. Makeenko. Its partition
function is also a KP hierarchy tau function, and this model is in turn equivalent to a Hermitian one-matrix model
with a general potential. Finally we prove that the generating function for general two-profile
Belyi morphisms is a GKM thus proving that it is also a KP hierarchy tau function in proper times.


\end{abstract}

\maketitle

\section{Introduction}\label{s:intro}
\setcounter{equation}{0}
In general, Hurwitz numbers pertain to combinatorial classes of ramified mappings $f:{\mathbb CP^1\to\Sigma_g}$
of the complex projective line onto a Riemann surface of genus $g$. Commonly, single and double Hurwitz numbers
correspond to the cases in which ramification profiles (defined by the corresponding Young tableauxes
$\lambda$ or $\lambda$ and $\mu$) are respectively given at one ($\infty$) or two ($\infty$ and $1$) distinct points 
whereas we assume the existence of $m$ other distinct ramification points with only simple ramifications.

Generating functions for Hurwitz numbers have been considered for long in mathematical physics. Notably, Okounkov and
Pandharipande \cite{OP} showed that the exponential of the generating function for
double Hurwitz numbers is a tau-function of the Kadomtsev--Petviashvili (KP) hierarchy. The same result was
obtained by A.~Yu.~Orlov and Shcherbin \cite{OS}, \cite{Or} using the Schur function technique and, in a more
general setting, by Goulden and Jackson \cite{GJ} using Plucker relations. 

Orlov and Shcherbin \cite{OS} also addressed the case of the generating function for the case of
Grothendieck {\em dessins d'enfants} where we have only three ramification points with multiple 
ramifications and the ramification profile is fixed at one or two of these points. In this case, they
also obtained that the exponentials of the corresponding generating functions are the tau functions
of the KP hierarchy.

On the other hand, Hurwitz numbers manifest properties intrinsic for conformal theories including sets of 
Virasoro constraints and closely related loop equations.
That simple Hurwitz numbers satisfy the topological recursion---the technique originated in matrix models---was 
conjectures in \cite{Mar} and proved in \cite{EMS}. \cite{MM}. 

In a nice recent paper \cite{Zograf} Zograf provided recursion relations for the generating function
of Grothendieck's {\it dessins d'enfants} enumerating the Belyi pairs $(C,f)$, where $C$ is a smooth algebraic curve
and $f$ a meromorphic function $f:C\to \mathbb CP^1$ ramified only over the points $0,1,\infty\in \mathbb CP^1$.

We recall some mathematical results relating Belyi pairs to Galois groups and begin with

\begin{theorem}\label{thm:Belyi}{\rm (Belyi, \cite{Belyi})}  A smooth complex algebraic curve $C$ is defined
over the field of algebraic numbers $\overline{\mathbb Q}$ if and only if it exists a nonconstant meromorphic function $f$ on $C$ $(f:C\to {\mathbb C}P^1)$ ramified only over the points $0,1,\infty\in {\mathbb C}P^1$.
\end{theorem}

For a Belyi pair $(C,f)$ let $g$ be the genus of $C$ and $d$ the degree of $f$. If we take the inverse image
$f^{-1}([0,1])\subset C$ of the real line segment $[0,1]\in {\mathbb C}P^1$ we obtain a connected bipartite
 fat graph with $d$ edges with vertices being preimages of $0$ and $1$ and with the cyclic ordering of edges entering
 a vertex coming from the orientation of the curve $C$. This led Grothendieck to formulating the following lemma:

\begin{lm}\label{lm:Grot} {\rm (Grothendieck, \cite{Grot})}
There is a one-to-one correspondence between the isomorphism classes of Belyi pairs and connected bipartite fat graphs.
\end{lm}

We define a Grothendieck {\it dessin d'enfant} to be a connected bipartite fat graph representing a Belyi pair.

\begin{figure}[tb]
{\psset{unit=0.7}
\begin{pspicture}(-3,-2.4)(3,2.4)
\psframe[linecolor=white, fillstyle=solid, fillcolor=yellow](-3.5,0)(3.5,-2.5)
\pcline[linestyle=solid, linewidth=1pt](-3.5,0)(3.5,0)
\pcline[linestyle=solid, linewidth=2pt](-1,0)(1,0)
\rput(-1,0){\pscircle*{.1}}
\rput(1,0){\pscircle*{.1}}
\psarc[linecolor=white, linestyle=solid, linewidth=10pt](-1,0){2}{60}{300}
\psarc[linecolor=red, linestyle=dashed, linewidth=1.5pt](-1,0){1.8}{60}{300}
\psarc[linecolor=blue, linestyle=solid, linewidth=1pt](-1,0){2.2}{60}{300}
\psarc[linecolor=white, linestyle=solid, linewidth=10pt](1,0){2}{-120}{120}
\psarc[linecolor=green, linestyle=dotted, linewidth=2pt](1,0){1.8}{-120}{120}
\psarc[linecolor=blue, linestyle=solid, linewidth=1pt](1,0){2.2}{-120}{120}
\pcline[linecolor=white, linestyle=solid, linewidth=10pt](0,1.73)(0,-1.73)
\pcline[linecolor=red, linestyle=dashed, linewidth=1.5pt](-0.2,1.73)(-0.2,-1.73)
\pcline[linecolor=green, linestyle=dotted, linewidth=2pt](0.2,1.73)(0.2,-1.73)
\rput(0,1.73){\pscircle[linecolor=black, fillstyle=solid, fillcolor=white]{.25}}
\rput(0,-1.73){\pscircle*{.25}}
\rput(-4,0){\makebox(0,0)[cc]{\hbox{{$\infty_-$}}}}
\rput(4,0){\makebox(0,0)[cc]{\hbox{{$\infty_+$}}}}
\rput(-1,.4){\makebox(0,0)[cc]{\hbox{{$0$}}}}
\rput(1,.4){\makebox(0,0)[cc]{\hbox{{$1$}}}}
\rput(0,2.3){\makebox(0,0)[cc]{\hbox{{\small$\Lambda$}}}}
\rput(0,-2.3){\makebox(0,0)[cc]{\hbox{{\small$\overline\Lambda$}}}}
\end{pspicture}
}
\caption{\small The Belyi graph $\Gamma_1$ corresponding to the Belyi pair $(\mathbb CP^1,\hbox{id})$; $\infty_{\pm}$
indicate directions of approaching the infinite point in $\mathbb CP^1$.
By $\Lambda$, $\overline \Lambda$ we indicate the insertions of the external field in the matrix-model formalism
of Sec.~\ref{s:model}. For example, this graph contributes the term $N^2\beta\gamma \tr(\Lambda\overline\Lambda)$.}
\label{fi:Belyi}
\end{figure}
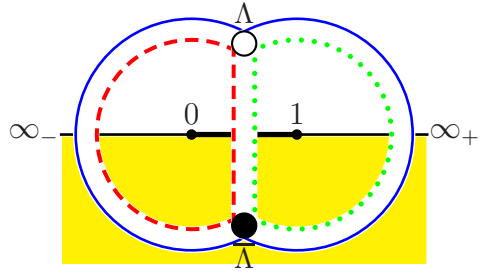

It is well known that we can naturally extend the dessin $f^{-1}([0,1])\subset C$ corresponding to a Belyi pair
$(C,f)$ to a bipartite triangulation of the curve $C$. For this, we cut the complex plane along the (real) line
containing $0,1,\infty$ coloring upper half plane white and lower half plane gray. This defines the partition of
$C$ into white and grey triangles such that white triangles has common edges only with grey triangles. We then consider a dual graph in which edges are of three types (pre-images of the three edges shown in Fig.~\ref{fi:Belyi}):
the type of an edge depend on which segment---$f^{-1}([0,1])\subset C$, $f^{-1}([1,\infty_+])\subset C$, or
$f^{-1}([\infty_-,0])\subset C$---it intersects ($\infty_{\pm}$ indicate the directions of approaching the
point of infinity along the real axis in $\mathbb CP^1$). 
Each face of the dual partition then contains a preimage of exactly
one of the points $0,1,\infty$, so they are of three sorts (bordered by solid, dotted, or dashed lines in the figure). We call such a graph a {\it Belyi fat graph}.

The type of ramification at infinity is determined by the set of solid-line bounded faces of a Belyi fat graph: the order of branching is $r$ for a $2r$-gon, so we introduce the generating function that
distinguishes between different types of
branching at infinity. We let $n_1,n_2,n_3$ denote the numbers of respective solid-, dotted-, and dashed-line
cycles (faces) and let $m_r$ denote the number of solid-line cycles of length $2r$ in a Belyi fat graph

We are interested in the following {\bf counting problem}: we are going to calculate the generating function
\be
{\mathcal F}\bigl[\{t_m\},\beta, \gamma;N\bigr]=\sum_{\Gamma}\frac{1}
{|\hbox{Aut\,}\Gamma|}N^{2-2g}\beta^{n_2}\gamma^{n_3}\prod_{i=1}^{n_1}t_{r_i},
\label{gen-fun}
\ee
where $N$, $\beta,\gamma$, and
$t_r$ are formal independent parameters and the sum ranges all (connected) Belyi fat graphs. Often a factor
$\alpha^{n_1}$ is also added; it can however be adsorbed into the times $t_r$ by scaling $t_r\to \alpha t_r$
for all $r$.

The structure of the paper is as follows.
In Sec.~\ref{s:model}, we show that generating function (\ref{gen-fun}) is the free energy of a special matrix model.
We demonstrate that this model is the two-logarithm matrix model of \cite{ChPal}, and it therefore
belongs to the class of generalized Kontsevich models (GKM) \cite{MMM}.
In Sec.~\ref{s:2-log}, we present
the solution of this model from paper \cite{ChPal} in which it was reduced, upon a special transformation of times,
to a Hermitian one-matrix model with a general potential. In Sec.~\ref{s:CEO}, we present the direct solution of
the original generating function in terms of the Hermitian one-matrix model without appealing to the external field
model thus again establishing the equivalence between the two models and describing the corresponding topological
recursion.
In Sec.~\ref{s:clean}, we construct the matrix model for
{\it clean} Belyi morphisms (those having ramifications only of type $(2,2,\dots,2)$ over $1$) and show that the corresponding generating function is the original Kontsevich--Penner model of \cite{ChM}. This model is also
equivalent \cite{ChM2} to the Hermitian one-matrix model with a general potential
and to the BGW model of \cite{MMS}. Finally, in Sec.~\ref{s:general},
we combine the techniques of Secs.~\ref{s:model},~\ref{s:2-log}, and~\ref{s:CEO} establishing that the
generating function for the
two-profile Belyi morphisms (with the given ramifications at two points, $\infty$ and $1$)
is again given by the GKM integral thus being a tau function of the KP hierarchy (that is, it
satisfies the bilinear Hirota relations). We conclude with the discussion of
our results.

Throughout the entire text we disregard all multipliers not depending on external fields; all equalities in the paper
must therefore be understood modulo such factors.

\section{The model}\label{s:model}
\setcounter{equation}{0}

In our conventions the indices $i$, $i_1$, $i_2$, etc. take positive integer values between $1$ and $\alpha N$,
the indices $j$, $j_1$, etc. take positive integer values between $1$ and $\beta N$, and the indices $k$, $k_1$, etc.
take positive integer values between $1$ and $\gamma N$. We introduce three complex-valued rectangular matrices
$R_{k,i}$, $G_{i,j}$, and $B_{j,k}$ and one diagonal matrix (the external field) $\Lambda_{i_1,i_2}=
\lambda_{i_1}\delta_{i_1,i_2}$. The action is given by the integral
\be
{\mathcal F}[\{t_r\},\beta,\gamma;N]:=
\int DR\,D\overline R\,DB\,D\overline B\,DG\,D\overline G\, e^{N\text{tr}(-B\overline B-R\overline R-G\overline G+
R\Lambda GB+\overline B\,\overline G\,\overline \Lambda\,\overline R)}.
\label{int1}
\ee
The free energy ${\mathcal F}[\{t_r\},\beta,\gamma;N]$
is given by the sum over all connected bipartite three-valent fat graphs $\Gamma$ weighted by
\be
\frac{1}{|\hbox{Aut\,}\Gamma|}N^{2-2g}\beta^{n_2}\gamma^{n_3}\prod_{r}t_r^{m_r}\qquad \Bigl(\sum_r m_r=n_1\Bigr)
\label{factor}
\ee
where $n_{1,2,3}$ are the respective numbers of solid-, dotted-, and dashed-line cycles in $\Gamma$,
\be
t_r:=\sum_{i=1}^{\alpha N}|\lambda_i|^{2r}
\label{t-r}
\ee
are the {\it times} of the model, and $m_r$ is the number of solid-line cycles of length $2r$ in $\Gamma$.
Measures of integration are the standard Haar measures; for instance,
$$
DR\,D\overline R:=\prod_{k=1}^{\gamma N}\prod_{i=1}^{\alpha N}d\hbox{Re}\,R_{k,i}d\hbox{Im}\,R_{k,i}.
$$
The logarithm of the integral (\ref{int1}) is therefore just the
generating function (\ref{gen-fun}) for the Belyi graphs.

Integrating w.r.t. $B,\overline B$ we obtain the integral
\be
\int DR\,D\overline R\,DG\,D\overline G\, e^{N\text{tr}(-R\overline R-G\overline G+
R\Lambda G\overline G\,\overline \Lambda\,\overline R)}
\label{int2}
\ee
in which we can perform the Gaussian integration w.r.t. $G,\overline G$ thus obtaining
\be
\int DR\,D\overline R\, e^{-N\text{tr}(R\overline R)}
\det\left[\delta_{i_1,i_2}-(\overline \Lambda\,\overline R R\Lambda)_{i_1,i_2}\right]^{-\beta N}.
\label{int3}
\ee
After the change of variables $R\to R\Lambda$ this integral becomes
\be
\prod_{i=1}^{\alpha N}|\lambda_i|^{-2\gamma N}
\int DR\,D\overline R\, e^{-N\text{tr}\bigl(\overline R R [\Lambda\overline\Lambda]^{-1}\bigr)}
\det\left[\delta_{i_1,i_2}-(\overline R R)_{i_1,i_2}\right]^{-\beta N}.
\label{int31}
\ee
For definiteness, let $\gamma\ge \alpha$. A general rectangular matrix $\overline R$ can then be reduced to
the form $\overline R=U^\dagger \overline M V$, where $U\in U(\alpha N)$, $V\in U(\gamma N)/U((\gamma-\alpha)N)$, and
$$
\overline M=\left(\begin{array}{ccc|cc}
          \overline m_1 & 0 & 0 & 0&0 \\
          0 & \ddots & 0 & 0&0 \\
           0 & 0 &  \overline m_{\alpha N}& 0&0
        \end{array}
\right).
$$
In the vicinity of the unities of the unitary groups, we can write $U=e^{i\epsilon H}$ and $V=e^{i\epsilon Q}$ with the  Hermitian $(\alpha N\times \alpha N)$-matrix $H$ and Hermitian $(\gamma N\times \gamma N)$-matrix $Q$ of the form
\be
Q=\left(
    \begin{array}{c|c}
      \widetilde H & P \\
      \hline
      P^\dagger & 0 \\
    \end{array}
  \right),
  \label{Q}
\ee
in which $\widetilde H$ is another Hermitian $(\alpha N\times \alpha N)$-matrix and $P$ is the general complex
$(\alpha N\times (\gamma-\alpha)N)$-matrix.
The Jacobian of the transformation
\be
D\overline R\,DR=\hbox{Jac}\, DU\,DV\,\prod_i dm_{i}d\overline m_{i}
\label{RR2UVM}
\ee
can then be easily
calculated (see Appendix~A) to be
\be
\hbox{Jac}=
\prod_{1\le i_1<i_2\le \alpha N}(|m_{i_2}|^2-|m_{i_1}|^2)^2\prod_{i=1}^{\alpha N}|m_i|^{2(\gamma-\alpha)N},
\label{Jac}
\ee
Introducing the new variables $x_i=|m_i|^2$ ranging from zero to infinity, we reduce the integral in (\ref{int31})
to the $\alpha N$-fold integral w.r.t. $x_i$ and to the integration w.r.t. the unitary group:
\bea
&{}&\prod_{i=1}^{\alpha N}|\lambda_i|^{-2\gamma N}
\int_{0}^\infty dx_1\dots dx_{\alpha N}\Bigl[\int DU e^{-N\sum_{i_1,i_2} x_{i_1}U_{i_1,i_2}|\lambda_{i_2}|^{-2}U^\dagger_{i_2,i_1}}\Bigr]\times\nonumber\\
&{}&\qquad\times
\bigl[\Delta(x)\bigr]^2 \prod_{i=1}^{\alpha N}\bigl[x_i^{(\gamma-\alpha)N}(1-x_i)^{-\beta N}\bigr].
\label{int4}
\eea
The integral over $DU$ is given by the Itzykson--Zuber--Mehta formula
(we write it having in mind that we subsequently integrate it
over variables $x_i$ with a totally symmetric measure),
$$
\int DU\, e^{-N\sum_{i_1,i_2} x_{i_1}U_{i_1,i_2}|\lambda_{i_2}|^{-2}U^\dagger_{i_2,i_1}}=
\frac{e^{-N\sum_i x_i |\lambda_i|^{-2}}}{\Delta(x_i)\Delta(|\lambda_i|^{-2})},
$$
so the final formula for the generating function reads
\be
\frac{\prod_{i=1}^{\alpha N}|\lambda_i|^{-2\gamma N}}{\Delta(|\lambda|^{-2})}
\int_{0}^\infty dx_1\dots dx_{\alpha N}\Delta(x)
e^{N\sum_{i} [-x_{i}|\lambda_{i}|^{-2}+(\gamma-\alpha)\log x_i -\beta\log (1-x_i)]}.
\label{int5}
\ee

The integral (\ref{int5}) is equivalent to the matrix-model integral
\be
\prod_{i=1}^{\alpha N}|\lambda_i|^{-2\gamma N}
\int_{\alpha N\times \alpha N} DH_{\ge 0}
e^{N\tr[-H\Lambda^{-2}+(\gamma-\alpha)\log H-\beta \log(1-H)]},
\label{int6}
\ee
where the integration goes over Hermitian $(\alpha N\times \alpha N)$-matrices with positive eigenvalues.
We thus obtain the following statement.

\begin{lm}\label{lm:gen-Groth}
The generating function for Grothendieck dessins d'enfants (Belyi fat graphs (\ref{gen-fun}))
is the matrix-model integral (\ref{int6}).
\end{lm}

The integral (\ref{int6})
belongs to the class of {\it generalized Kontsevich models} (GKM) \cite{MMM}; in terms of variables
$\xi_i=1/|\lambda_i|^2$ it can be calculated as the ratio of determinants of
$(\alpha N\times \alpha N)$-matrices,
$$
\Bigl\| \frac{\partial^{i_1-1}f(\xi_{i_2})}{\partial \xi_{i_2}^{i_1-1}}\Bigr\|/\Delta(\xi),
$$
where
$$
f(\xi)=\int_0^{\infty} dx e^{-Nx\xi}x^{(\gamma-\alpha)N}(1-x)^{-\beta N},
$$
and as such is a tau-function of the Kadomtsev--Petviashvili (KP) hierarchy in times
$t_n=\sum_i \xi_i^{-n}=\sum_i |\lambda_i|^{2n}$ (cf. (\ref{t-r}))
i.e., we come to the following theorem
proved by Zograf \cite{Zograf}  by purely combinatorial means with the using of the cut-and-joint operator.

\begin{theorem}\label{thm:Zograf}
The generating function for Belyi fat graphs (\ref{gen-fun}) is the tau-function of the
KP hierarchy in times (\ref{t-r}).
\end{theorem}

The integral (\ref{int6}) was studied by one of the authors and Palamarchuk \cite{ChPal} in relation to exploring possible explicit solutions of matrix models with external fields. It was called the {\it two-logarithm model} there and it was
proved that this integral admits Virasoro constraints that, upon a proper change of times, become the
Virasoro constraints of the matrix model introduced in \cite{ChM}
(the term Kontsevich--Penner model was coined there), which, in turn, is equivalent~\cite{ChM2}
to a Hermitian one-matrix model with the potential related to the external-field variables
$\xi_i$ via the Miwa transformation. As such, this integral must also satisfy the equations of the Toda chain hierarchy.

\begin{remark}\label{rm:no-limit}
An important remark concerning integral (\ref{int6}) is that its asymptotic behavior as $N\to\infty$ is different
depending on whether $\gamma-\alpha\simeq O(1)$ or $\gamma-\alpha\simeq O(1/N)$.
In the first case, we have an infinite repulsive
potential at the origin and an eigenvalue distribution is confined within an interval $[x'_-,x'_+]$ (see below)
with $0<x'_-<x'_+$. The $1/N$-expansion then is ``insensitive''
to the hard edge at the origin, and we can assume that we
integrate over the whole real axis (the difference between the restricted and nonrestricted integrations is then
exponentially small in $N$). If $\gamma=\alpha$ or $\gamma-\alpha\sim O(1/N)$,
representation (\ref{int6}) still remains valid, but in this case
the eigenvalue support is $[0,x'_+]$, so it reaches the hard edge $x=0$ at the origin. We then again have a topological expansion (about $1/N$-expansion in matrix models with hard edges,
see, e.g., review \cite{Ch06}) but with the differential $ydx$ finite at $x=0$
($y\sim1/\sqrt{x}$ as $x\to 0$ and $y\sim \sqrt{x-x'_+}$ as $x\to x'_+$). The asymptotic expansions of
 integral (\ref{int6}) are therefore different in the corresponding regimes and do not admit an analytical transition
 as $\gamma\to\alpha$.
\end{remark}

\begin{remark}\label{rm:explanation}
In Sec.~\ref{s:CEO}, we present a simpler, straightforward way of proving that generating function
(\ref{gen-fun}) for general Belyi morphisms is indeed a Hermitian one-matrix model free energy. However, the
external field technique of this and next sections will be instrumental when proving a general correspondence
between the generating functions for clean (Sec.~\ref{s:clean}) and two-profile (Sec.~\ref{s:general})
Belyi morphisms and free energies of the corresponding generalized Kontsevich models.
\end{remark}

\section{The two-logarithm matrix model}\label{s:2-log}
\setcounter{equation}{0}
In this section, we present the results of \cite{ChPal} adapted to the notation of integral (\ref{int6}).

\subsection{Constraint equations for integral (\ref{int6})}
We first perform the variable change
\be
\begin{array}{lll}
\tilde N=\alpha N, & \tilde \Lambda=\Lambda^{-2}/(2\alpha), & \tilde H=2H-1\\
\tilde\alpha=\beta/\alpha, & \tilde \beta=1-\gamma/\alpha. &
\end{array}
\label{relations}
\ee
in (\ref{int6}). Disregarding here and hereafter
factors not depending on $\lambda$'s, the integral then takes the form
\be
\prod_{i=1}^{\tilde N}\Bigl[|\tilde\lambda_i|^{\gamma N}e^{-\tilde N|\tilde\lambda_i|}\Bigr]
\int_{\tilde N\times \tilde N} D\tilde H_{\ge 0}
e^{-\tilde N\tr [\tilde H\tilde \Lambda +
\tilde \alpha \log(1-\tilde H)+\tilde \beta\log(1+\tilde H)]}
:=\prod_{i=1}^{\tilde N}\Bigl[|\tilde\lambda_i|^{\gamma N}e^{-\tilde N|\tilde\lambda_i|}\Bigr]
{\mathcal Z}[\tilde\lambda],
\label{int7}
\ee
where we let ${\mathcal Z}[\tilde\lambda]$ denote the integral (\ref{int6}) without the normalization factor.

The Schwinger--Dyson equations for the integral (\ref{int7}) follow from the identity (here all the indices range
from $1$ to $\alpha N$)
\be
\left(\frac1{{\tilde N}^3}\frac{\partial}{\partial \tilde\Lambda_{jk}}
\frac{\partial}{\partial \tilde\Lambda_{li}}-\frac{1}{\tilde N}\right)
\int_{\tilde N\times \tilde N} D\tilde H\,\frac{\partial}{\partial \tilde H_{ij}}e^{-\tilde N\tr [\tilde H\tilde \Lambda +
\tilde \alpha \log(1-\tilde H)+\tilde \beta\log(1+\tilde H)]}=0.
\label{SD-eq}
\ee
In terms of the eigenvalues $\tilde\lambda_i$ of the matrix $\tilde \Lambda$, the
corresponding $\tilde N$ equations read
\be
\left[-\frac{1}{\tilde N^2}\tilde\lambda_{i_1}\frac{\partial^2}{\partial\tilde\lambda_{i_1}^2}-
\frac{1}{\tilde N^2}\sum_{i_2\ne i_1}\frac{\tilde\lambda_{i_2}}{\tilde\lambda_{i_2}-\tilde\lambda_{i_1}}
\Bigl(\frac{\partial}{\partial\tilde\lambda_{i_2}}-\frac{\partial}{\partial\tilde\lambda_{i_1}}\Bigr)+
\frac{\tilde\alpha+\tilde\beta-2}{\tilde N}\frac{\partial}{\partial\tilde\lambda_{i_1}}+\tilde\beta-\tilde\alpha+
\tilde\lambda_{i_1}
\right]{\mathcal Z}[\tilde\lambda]=0,
\label{SD-lambda}
\ee

We can equivalently write the constraint equations (\ref{SD-lambda}) in terms of the {\it times}
\be
t_n=\frac 1n \sum_i \frac1{{\tilde \lambda_i}^n}, \quad n\ge 1.
\label{lambda-times}
\ee
They then becomes the set of {\it Virasoro constraints}\footnote{The authors were reported by M.~Kazarian that the
same constraints can be derived by pure combinatorial means [M.~Kazarian, P.~Zograf, paper in preparation].}
\be
V_k{\mathcal Z}\bigl(\{t_n\}\bigr)=0,\quad k\ge 0,
\ee
where
\bea
V_k[t]&:=&-\sum_{m=1}^\infty mt_m\frac{\partial}{\partial t_{m+k}}-\sum_{m=1}^k \frac{\partial}{\partial t_{m}}
\frac{\partial}{\partial t_{k-m}}-\tilde N (\tilde \alpha-\tilde\beta+1)(1-\delta_{k,0}-\delta_{k,-1})
\frac{\partial}{\partial t_{k}}\nonumber\\
&{}&+\bigl[2\tilde N(1-\delta_{k,-1})+\delta_{k,-1} t_1\bigr]\frac{\partial}{\partial t_{k+1}}
+\tilde N^2\tilde\alpha(\tilde \beta-1)\delta_{k,0},\quad k=-1,0,1,\dots.
\label{Vk}
\eea
(Here, for the future use, we have also introduced the operator $V_{-1}$.)

The operators $V_k$ enjoy the Virasoro algebra
\be
[V_k,V_l]=(l-k)V_{k+l},\quad k,l\ge -1.
\ee

\subsection{Equivalence to the Hermitian one-matrix model}\label{ss:equivalence}

In \cite{ChPal} it was shown that the two-logarithm model is equivalent to the Kontsevich--Penner model \cite{ChM},
which in turn was known \cite{ChM2}, \cite{MMM} to be equivalent to a Hermitian one-matrix model. In this paper,
we skip the intermediate step and demonstrate the equivalence between (\ref{int6}) and a Hermitian one-matrix
model defined as an integral
\be
{\mathcal Z}_{\text{1MM}}\bigl[\{\xi_m\},M\bigr]:=\int_{M\times M}DY\,e^{-V(Y)},\quad V(Y)=\sum_{m=1}^\infty \xi_m\tr Y^m.
\label{1MM}
\ee
It is well-known that this integral satisfies the set of Virasoro constraints uniformly written in the form
\be
L_n {\mathcal Z}_{\text{1MM}}\bigl[\{\xi_m\},M\bigr]=\Bigl\{
\sum_{m=0}^n\frac{\partial^2}{\partial \xi_m\partial\xi_{n-m}}
+\sum_{m=1}^\infty m\xi_m\frac{\partial}{\partial \xi_{n+m}} \Bigr\}
{\mathcal Z}_{\text{1MM}}\bigl[\{\xi_m\},M\bigr]=0,
\ n\ge -1,
\label{Virasoro-1MM}
\ee
where we have used a convenient notation
$\frac{\partial}{\partial \xi_{0}}{\mathcal Z}_{\text{1MM}}\bigl[\{\xi_m\},M\bigr]
=-M {\mathcal Z}_{\text{1MM}}\bigl[\{\xi_m\},M\bigr]$.

In order to establish the correspondence it is necessary to shift the original variable $\tilde\lambda$,
\be
\mu_i=\tilde\lambda_i-\rho,\quad \rho\in{\mathbb C},
\label{shift}
\ee
introducing an auxiliary parameter $\rho$. We also introduce the new times
\be
\tau_n:=\frac1n \sum_{i=1}^{\tilde N}\frac1{\mu_i^n},\quad n\ge1,
\label{tau}
\ee
and the new normalizing factor
\be
{\mathcal N}[\mu]:=\prod_{i=1}^{\tilde N}\left[\mu_i^{\tilde N(\tilde \beta-1)}e^{\tilde N\mu_i}\right]
\label{N}
\ee
The following set of constraints was found in \cite{ChPal}:
\begin{lm}\label{lm:Lk}{\rm(see \cite{ChPal})}
The normalized integral ${\mathcal Z}[\tilde\lambda]/{\mathcal N}[\mu]$ where $\tilde\lambda_i=\mu_i+\rho$
satisfies the set of Virasoro constraints
$$
{\mathcal L}_k \bigl[{\mathcal Z}[\tilde\lambda]/{\mathcal N}[\mu]\bigr]=0,\ k=-1,0,1,\dots,
$$
in times (\ref{tau}) with
\bea
{\mathcal L}_k&=&-\sum_{m=1+\delta_{k,-1}}^\infty m(\tau_m-2\tilde N\delta_{m,1})
\frac{\partial}{\partial\tau_{m+k}}
-\sum_{m=1}^{k-1}\frac{\partial^2}{\partial \tau_m\partial\tau_{k-m}}
+2\tilde N\alpha_{\text{KP}}(1-\delta_{k,0}-\delta_{k,1})\frac{\partial}{\partial \tau_{k}}\nonumber\\
&{}&-2\varphi \tilde N\sum_{m=1+\delta_{k,-1}}\frac{1}{(-\rho)^m}\frac{\partial}{\partial \tau_{k+m}}
-(\tilde N\alpha_{\text{KP}})^2\delta_{k,0}
+\tilde N\alpha_{\text{KP}}\Bigl(\tau_1-2\tilde N-\frac{2\varphi\tilde N}{\rho}\Bigr)\delta_{k,-1},
\label{Lk}
\eea
where $\alpha_{\text{KP}}=\tilde\beta-1$ and $\varphi=-(\tilde\alpha+\tilde\beta-1)/2$.
\end{lm}

\begin{remark}\label{rm:Virasoro}
In order to derive constraints (\ref{Lk}) the following trick was used in \cite{ChPal}: constraint equations
(\ref{SD-lambda}) after shift (\ref{shift}) were written in the form
$$
\sum_{k=1}^\infty \mu_i^{-k} L_k{\mathcal Z}[\tilde\lambda]=0,
$$
where
\bea
L_k&=& V_{k+1}[\tau] +\rho V_k[\tau] +\rho\tilde N (\tilde\alpha+\tilde\beta-1)
\Bigl( (1-\delta_{k,0}-\delta_{k,-1})\frac{\partial}{\partial\tau_k}-(\tilde\beta-1)\tilde N\delta_{k,0} \Bigr)
\nonumber\\
&{}&+\rho(\tilde\beta-1)\tilde N (\tau_1-2\tilde N) \delta_{k,-1},\quad k\ge -1,\nonumber
\eea
were differential operators in (shifted) times $\tau_s$ and where we let
$V_s[\tau]$ denote operators (\ref{Vk}) upon the
substitution $t\to \tau$. The ``proper'' Virasoro operators ${\mathcal L}_k$ (\ref{Lk}) were finally obtained
upon the upper-triangular transformation
$$
{\mathcal L}_k=\sum_{s=0}^\infty \frac{(-1)^s}{\rho^{s+1}}L_{k+s},\quad k\ge -1.
$$
We see that in order to perform all these replacements we have to keep $\rho$ nonzero and finite.
\end{remark}

\begin{lm}\label{lm:sub}{\rm(see \cite{ChPal})}
Upon the substitution
\be
\xi_n=\tau_n+\frac 1n \frac{2\varphi\tilde N}{(-\rho)^n}-2\tilde N\delta_{n,1},\quad M=\tilde N\alpha_{\text{KP}}
\label{time-change}
\ee
the Virasoro constraints (\ref{Lk}) become the Virasoro constraints (\ref{Virasoro-1MM}) of the Hermitian
one-matrix model. Because these conditions determine the corresponding integrals unambiguously, these two
models are equivalent.
\end{lm}

In terms of the original variables, we have the following lemma.

\begin{lm}\label{lm:equivalence}
The generating function ${\mathcal F}[\{t_r\},\beta,\gamma;N]$ (\ref{gen-fun}) for the Belyi fat graphs
is given by the exact formula
\bea
e^{{\mathcal F}[\{t_r\},\beta,\gamma;N]}&=&\prod_{i=1}^{\alpha N}
\left[\Bigl(\frac1{2\alpha}-\rho|\lambda_i|^2\Bigr)^{-\gamma N}
e^{\alpha N\Bigl(\frac{1}{2\alpha |\lambda_i|^2}-\rho\Bigr)}\right]\times\nonumber\\
&{}&\times {\mathcal Z}_{\text{1MM}}
\Bigl[\xi_m=\tau_m+\frac 1m \frac{(\gamma-\beta)N}{(-\rho)^m}-2\alpha N\delta_{n,1},\,M=-\gamma N\Bigr]
\label{F=1MM}
\eea
with $\tau_m=\frac 1m \sum_{j=1}^{\alpha N}\frac{1}{\mu_j^m}$ where $\mu_i+\rho=1/(2\alpha |\lambda_i|^2)$.
Here ${\mathcal Z}_{\text{1MM}}\bigl[\{\xi_m\},M\bigr]$ is matrix integral (\ref{1MM}).
\end{lm}

In the next section we demonstrate that this statement enables us to write explicit formulas for terms of
the genus expansion of ${\mathcal F}$ provided we know the answer for the free energy of matrix model (\ref{1MM})
either in terms of momentums \cite{ACKM} or in terms of the topological recursion technique of
\cite{Ey}, \cite{ChEy}, \cite{CEO}, \cite{AlMM}.

\begin{remark}
The shift of variables (\ref{shift}) is a convenient technical tool that was used in \cite{ChPal} for passing to the
full half-Virasoro constraint algebra that includes also
the operator $L_{-1}$. If $|\gamma-\alpha|\lesssim O(1/N)$ we have a hard edge at the origin, which is 
specific for the complex matrix model of \cite{AKM} or the BGW model of \cite{MMS}, and we shall lose the
$L_{-1}$ Virasoro operator.\footnote{The authors thank A.~Mironov for this comment.} We reconstruct the 
$L_{-1}$-operator in the model with logarithmic potential for the price of unfreezing all times of the hierarchy.
And, as we demonstrate in the next section, the
final answers for genus expansion terms do not depend on the auxiliary parameter $\rho$.
\end{remark}

\subsection{The genus expansion}\label{ss:large-N}

An extensive literature is devoted to solving the one-matrix model (\ref{1MM}) in the topological (genus)
expansion; its free energy $F$ admits a representation $F=\sum_{h=0}^\infty M^{2-2h}F_h$, which can be
interpreted as a semiclassical expansion of a (quasi)stationary statistical theory. As such, in the large-$M$ limit,
we observe a stationary distribution of eigenvalues described by a {\it spectral curve} of the model. In the present
paper, as in \cite{ChPal}, we assume that this stationary distribution spans a single interval, and we therefore have
a {\it one-cut} solution based on a spectral curve that is just a double cover of the complex plane with
two branching points, $x_+$ and $x_-$ (a sphere).
These two points are determined by the constraint equations for the so-called
master loop equation \cite{Migdal}
\be
\oint_{C_D}\frac{dw}{2\pi i}\frac{V'(w)}{\sqrt{(w-x_+)(w-x_-)}}=0,\qquad
\oint_{C_D}\frac{dw}{2\pi i}\frac{wV'(w)}{\sqrt{(w-x_+)(w-x_-)}}=2M,
\label{loop}
\ee
where the integration contour encircles the eigenvalue domain (the interval $[x_-,x_+]$ in this case) and not other
singularities (including possible singularities of $V'(w)$).

After the Miwa time transformation (\ref{time-change}) we obtain for $V'(w)$ the expression
\be
V'(w)=-2\alpha N-\sum_{i=1}^{\alpha N}\frac{1}{w-\mu_i}-(\gamma-\beta)N\frac {1}{w+\rho}
\label{Vprime}
\ee
and we assume that all $\mu_i$ and $-\rho$ are situated
outside the integration contour. We can then take the integrals in
(\ref{loop}) by residues at $\mu_i$, $-\rho$, and infinity. For the first equation we obtain
$$
-2\alpha N+\sum_{i=1}^{\alpha N}\frac{1}{\sqrt{(\mu_i-x_+)(\mu_i-x_-)}}
+(\gamma-\beta)N\frac{1}{\sqrt{(p+x_+)(p+x_-)}}=0
$$
and shifting the branching points
$$
x_++\rho=x'_+,\quad x_-+\rho=x'_-
$$
and recalling that $\mu_i+\rho=\tilde \lambda_i$ we obtain the constraint equation solely in terms
of $\tilde\lambda_i$:
\be
-2\alpha N+\sum_{i=1}^{\alpha N}\frac{1}{\sqrt{(\tilde\lambda_i-x'_+)(\tilde\lambda_i-x'_-)}}
+(\gamma-\beta)N\frac{1}{\sqrt{x'_+x'_-}}=0
\label{constraint1}
\ee
For the second constraint equation we obtain
$$
-\alpha N(x'_++x'_--2\rho)
+\sum_{i=1}^{\alpha N}\frac{\tilde\lambda_i-\rho}{\sqrt{(\tilde\lambda_i-x'_+)(\tilde\lambda_i-x'_-)}}-\alpha N
-(\gamma-\beta)N +(\gamma-\beta)N\frac{-\rho}{\sqrt{x'_+x'_-}}=-2\gamma N
$$
and the term linear in $\rho$ is just the first constraint equation and thus vanishes. So, the second
constraint equation becomes
\be
(\gamma+\beta-\alpha)N-\alpha N(x'_++x'_-)
+\sum_{i=1}^{\alpha N}\frac{\tilde\lambda_i}{\sqrt{(\tilde\lambda_i-x'_+)(\tilde\lambda_i-x'_-)}}=0.
\label{constraint2}
\ee
We see that, as expected, all the dependence on $\rho$ disappears from constraint equations (\ref{constraint1})
and (\ref{constraint2}).

\begin{remark}\label{rm:correspondence}
Equations (\ref{constraint1}) and (\ref{constraint2}) exactly coincide with the respective first and second
constraint equations
in Eq.~(2.14) of \cite{ChPal} upon the substitution
\be
\begin{array}{lll}
\lambda\to\tilde\lambda, & N\to \alpha N & \beta-\alpha\to 1-\gamma/\alpha-\beta/\alpha,\\
c\to (\beta-\gamma)^2/4\alpha^2, & b/a\to -x'_+-x'_-, & c/a\to x'_-x'_+.
\end{array}
\label{correspondence}
\ee
The answer for ${\mathcal F}_0$ (formula (2.16) in \cite{ChPal})
obtained from these constraint equations therefore coincides (up to the normalization factor
$\prod_{i=1}^{\alpha N}\Bigl[|\tilde\lambda_i|^{\gamma N}e^{-\alpha N|\tilde\lambda_i|}\Bigr]$) with the genus zero contribution to generating function (\ref{gen-fun}).
\end{remark}

\subsubsection{Genus-zero term}
It follows from Remark~\ref{rm:correspondence} that the genus-zero term ${\mathcal F}_0$ of our generating
function (\ref{gen-fun}) upon the substitutions (\ref{correspondence}) and (\ref{relations})
coincides with $F_0$ found in \cite{ChPal} with the added
normalization term $\sum_{i=1}^{\alpha N}\bigl[\gamma N\log \tilde\lambda_i-\alpha N \tilde\lambda_i\bigr]$.
In terms of variables $x'_\pm,\tilde\lambda$ the corresponding expression reads
\bea
{\mathcal F}_0&=&\frac14(\beta^2N^2+\gamma^2N^2)\log\bigl[(x'_+-x'_-)^2\bigr]\nonumber\\
&{}&
+N^2(\alpha-\beta-\gamma)\biggl[|\beta-\gamma|
\log\biggl(\frac{x'_++x'_--2\sqrt{x'_+x'_-}}{x'_++x'_-+2\sqrt{x'_+x'_-}}\biggr)
+\frac{x'_+ +x'_-}{2}\biggr]\nonumber\\
&{}&+N^2\Bigl[\frac{\alpha^2}{8}(x'_+ +x'_-)^2+\alpha|\beta-\gamma|\sqrt{x'_+x'_-}-
\frac{(\beta-\gamma)^2}{4}\log[x'_+x'_-]\Bigr]\nonumber\\
&{}&+N\sum_{i=1}^{\alpha N}\biggl\{\frac{\beta+\gamma}{2}\log|\tilde\lambda_i|+g(\tilde\lambda_i)-\tilde\lambda_i
+\frac{\alpha-\beta-\gamma}{2}\log\Bigl(\tilde\lambda_i-\frac{x'_+ +x'_-}{2}+g(\tilde\lambda_i)\Bigr)\biggr.
\nonumber\\
&{}&-\biggl. \frac{|\beta-\gamma|}{4}
\log\frac{g(\tilde\lambda_i)-\frac{\tilde\lambda_i(x'_++x'_-)}{2\sqrt{x'_+x'_-}}+\sqrt{x'_+x'_-}}
{g(\tilde\lambda_i)+\frac{\tilde\lambda_i(x'_++x'_-)}{2\sqrt{x'_+x'_-}}-\sqrt{x'_+x'_-}}
\biggr\}\nonumber\\
&{}&-\frac14\sum_{i_1,i_2=1}^{\alpha N}\log\Bigl[g(\tilde\lambda_{i_1})g(\tilde\lambda_{i_2})
+\tilde\lambda_{i_1} \tilde\lambda_{i_2}
-\frac{\tilde\lambda_{i_1}+\tilde\lambda_{i_2}}{2}(x'_+ +x'_-)+x'_+x'_-\Bigr]
\label{F0}
\eea
where we have introduced the notation $g(\tilde\lambda_i):=\sqrt{(\tilde\lambda_i-x'_+)(\tilde\lambda_i-x'_-)}$.

It is easy to see that in the domain of large $\tilde\lambda_i$, the expansion in (\ref{F0}) contains only negative
powers of $\tilde\lambda$: the linear and the logarithmic in $\tilde\lambda_i$ terms vanish in this domain.

\subsubsection{Higher genus expressions}
All higher genus corrections to the Hermitian one-matrix model can be written in terms of {\it moments}
\cite{ACKM} $M_r$, $J_r$ of the potential:
\be
M_r=\oint_{C_D}\frac{dw}{2\pi i}\frac{V'(w)}{(w-x_+)^{r+1/2}(w-x_-)^{1/2}},\qquad
J_r=\oint_{C_D}\frac{dw}{2\pi i}\frac{V'(w)}{(w-x_+)^{1/2}(w-x_-)^{r+1/2}}, \ r\ge 1.
\label{moments1}
\ee
Using representation (\ref{Vprime}), we obtain for the moments the following expressions
\be
\begin{array}{l}
M_r=\sum\limits_{i=1}^{\alpha N}\dfrac{1}{(\tilde\lambda_i-x'_+)^{r+1/2}(\tilde\lambda_i-x'_-)^{1/2}}
+(\gamma-\beta)N\dfrac{(-1)^r}{(x'_+)^{r+1/2}(x'_-)^{1/2}}\phantom{\Biggr(}\cr
J_r=\sum\limits_{i=1}^{\alpha N}\dfrac{1}{(\tilde\lambda_i-x'_+)^{1/2}(\tilde\lambda_i-x'_-)^{r+1/2}}
+(\gamma-\beta)N\dfrac{(-1)^r}{(x'_+)^{1/2}(x'_-)^{r+1/2}}\phantom{\Biggr(}
\end{array}
\quad r\ge 1.
\label{moments2}
\ee
After substitution (\ref{moments2}), the answer for ${\mathcal F}_h$ for generating function (\ref{gen-fun})
is given by that of the standard Hermitian one-matrix model. We have thus proved the following lemma

\begin{lm}\label{Fh}
In terms of moments (\ref{moments2}), every term ${\mathcal F}_h$ corresponding to the genus $h>0$
has a polynomial form for higher $h$ \cite{ACKM},
\be
{\mathcal F}_h=\sum_{r_s>1,q_s>1}\bigl\langle r_1\dots r_m;q_1\dots q_l|r\,q\,p\bigr\rangle_h
\frac{M_{r_1}\cdots M_{r_m}J_{q_1}\cdots J_{q_l}}{M_1^r J_1^q |x'_+-x'_-|^p},\quad h>1,
\label{F2}
\ee
and~\cite{ACM}
\be
{\mathcal F}_1=-\frac1{24}\log\bigl[M_1J_1|x'_+-x'_-|^4\bigr].
\label{F1}
\ee
Here $\bigl\langle r_1\dots r_m;q_1\dots q_l|r\,q\,p\bigr\rangle_h$ are finite (for a fixed $h$) sets of
rational numbers given by the topological recursion technique for the standard Hermitian one-matrix model
(see \cite{ChEy}). They are subject to restrictions: $m+l-r-q=2-2h$,
$\sum_{s=1}^m(r_1-1)+\sum_{s=1}^l(q_s-1)+p=4h-4$, $p\ge h-1$.
\end{lm}

Using topological recursion we can effectively calculate the numbers $\bigl\langle r_1\dots r_m;q_1\dots q_l|r\,q\,p\bigr\rangle_h$.
The quantity $|x'_+-x'_-|$, which is often denoted
by $d$, is the length of the interval of eigenvalue support.
Formulas (\ref{F2}), (\ref{F1}), and (\ref{moments2})
thus describe generating function (\ref{gen-fun}) in all orders of the genus expansion.

\section{Spectral curve and topological recursion}\label{s:CEO}
\setcounter{equation}{0}
In this section, we directly derive the spectral curve without appealing to a matrix model with external fields.
For this, we shrink all solid-line cycles assigning just the original times $t_r$ to the obtained $2r$-valent vertices of the field $B$, $\overline B$. The generating function (\ref{gen-fun}) is then described by the matrix-model integral over rectangular $(\gamma N\times \beta N)$-matrices $B$:
\be
{\mathcal Z}[t]=\int_{\gamma N\times \beta N}DB\,D\overline B\,
e^{-N\tr[B\overline B]+N\sum_{r=1}^\infty \frac 1r t_r\tr [(B\overline B)^r]},
\label{new1}
\ee
which, using the Jacobian from Appendix~A under assumption that $\beta>\gamma$, can be reduced to the
$\gamma N$-fold integral over positive $x_k$:
\be
{\mathcal Z}[t]=\int_0^\infty dx_1\dots dx_{\gamma N}[\Delta(x)]^2\prod_{k=1}^{\gamma N}x_k^{(\beta-\gamma)N}
e^{-N\sum\limits_{r=1}^\infty\sum\limits_{k=1}^{\gamma N} \frac 1r(\delta_{r,1}-t_r)x_k^r }.
\ee
This integral is again a Hermitian one-matrix model with a logarithmic term in the potential:
\be
{\mathcal Z}[t]=\int\limits_{\gamma N\times \gamma N} DX_{\ge 0}
e^{-N\tr\Bigl[\sum\limits_{r=1}^\infty \frac 1r(\delta_{r,1}-t_r)X^r-(\beta-\gamma)\log X\Bigr]},
\label{new2}
\ee
We have thus obtained another representation of generating function (\ref{gen-fun}).

\begin{lm}\label{lm:1MM}
Generating function (\ref{gen-fun}) can be presented as a Hermitian one-matrix model integral (\ref{new2}) with
a logarithmic term in the potential.
\end{lm}

Because we have reduced the original problem to a mere Hermitian one-matrix model integral, we can
directly apply a standard topological recursion
procedure \cite{ChEy} (see \cite{Ch} where it was generalized to the case of rational
functions $V'(x)$). We let
\be
U'(x):=N\sum\limits_{r=1}^\infty (\delta_{r,1}-t_r)x^{r-1}
\label{U}
\ee
denote the polynomial part of the potential with times $t_r$ with the shifted first time. The
hyperelliptic spectral
curve is a sphere with two branching points $x'_+$ and $x'_-$ whose positions are determined by
the standard constraints (\ref{loop}) in which
\be
V'(x)=U'(x)-\frac{N(\beta-\gamma)}{x},\quad M=\gamma N.
\ee
Constraints (\ref{loop}) then become
\be
\oint\limits_{C_D}\frac{dw}{2\pi i}\frac{U'(w)}{\sqrt{(w-x'_+)(w-x'_-)}}
=\frac{N(\beta-\gamma)}{\sqrt{x'_+x'_-}},\quad
\oint\limits_{C_D}\frac{dw}{2\pi i}\frac{wU'(w)}{\sqrt{(w-x'_+)(w-x'_-)}}=N(\beta+\gamma),
\label{loop-U}
\ee
i.e., precisely constraints (\ref{constraint1}) and (\ref{constraint2}) after the inverse Miwa
transformation.\footnote{The term $(\beta+\gamma)$ in the r.h.s. of the second equation is not a misprint.}

The $y$-variable of the topological recursion is given by the integral
over the contour that encircles the eigenvalue support and the point $x$,
\be
y(x):=\oint\limits_{C_{[x'_-,x'_+]}\cup\{x\} }\frac{dw}{2\pi i}
\frac{V'(w)\sqrt{(x-x'_+)(x-x'_-)}}{(w-x)\sqrt{(w-x'_+)(w-x'_-)}},
\ee
which can be evaluated by residues at infinity and at $w=0$ (due to the presence of a pole term in $V'(w)$) The
result reads
\be
y(x)=\biggl(\hbox{res}{}_\infty\Bigl[\frac{U'(w)}{(w-x)\sqrt{(w-x'_+)(w-x'_-)}}\Bigr]+
\frac{N(\beta-\gamma)}{\sqrt{x'_+x'_-}}\biggr)\sqrt{(x-x'_+)(x-x'_-)}
\ee

The {\bf genus expansion} for $h\ge1$ has the same form as in Lemma (\ref{Fh}) with the moments given by the
standard integrals taken by residues at infinity and at $w=0$:
\be
\begin{array}{l}
M_r=\hbox{res\,}_{w=\infty}\Bigl[\dfrac{U'(w)}{(w-x'_+)^{r+1/2}(w-x'_-)^{1/2}}\Bigr]
+(\gamma-\beta)N\dfrac{(-1)^r}{(x'_+)^{r+1/2}(x'_-)^{1/2}}\phantom{\Biggr(}\cr
J_r=\hbox{res\,}_{w=\infty}\Bigl[\dfrac{U'(w)}{(w-x'_+)^{1/2}(w-x'_-)^{r+1/2}}\Bigr]
+(\gamma-\beta)N\dfrac{(-1)^r}{(x'_+)^{1/2}(x'_-)^{r+1/2}}\phantom{\Biggr(}
\end{array}
\quad r\ge 1.
\label{moments2-1MM}
\ee
The term ${\mathcal F}_0$ has the general form \cite{BIPZ}
(for the number of eigenvalues equal $t_0 N$)
\be
{\mathcal F}_0=-\frac12 \int\limits_{C_{[x'_-,x'_+]}}y(x)V(x)-\zeta t_0,
\label{F0-1MM}
\ee
where $\zeta$ is the Lagrange multiplier most conveniently obtained as the limit of the integral
\be
\zeta=\lim_{\Lambda\to+\infty}\Bigl( \int_{x'_+}^{\Lambda} y(x)dx-V(\Lambda)-t_0\log\Lambda\Bigr).
\ee

\section{Generating functional for clean Belyi morphisms}\label{s:clean}
\setcounter{equation}{0}

\subsection{The model}\label{ss:model-clean}
A {\it clean} Belyi morphism is a special class of Belyi pairs $(C,f)$ that have profile
$(2,2,\dots,2,1,1,\dots,1)$ over the
branch point $1\in {\mathbb C}P^1$. This means that all dotted cycles (in Fig.~\ref{fi:Belyi}) have either lengths 2
(no ramification) or 4 (simple ramification). In \cite{DMSA} the authors demonstrated that the
generating function for ramifications of sort $(2,2,\dots,2)$
satisfies the topological recursion relations with the spectral curve $(x=z+z^{-1};\,y=z)$.

In this section, we demonstrate that the matrix model corresponding to clean Belyi morphisms is just
the Kontsevich--Penner model~\cite{ChM}, which is in turn equivalent~\cite{ChM2} to the Hermitian one-matrix model
with a general potential.

We thus have to calculate the generating function (\ref{gen-fun}) in which the sum ranges over
only clean Belyi morphisms.
In terms of the diagrammatic technique of Sec.~\ref{s:model} this means that we count only dotted cycles of lengths 2 and 4. Counting cycles of length 2 reduces to a mere changing of the normalization of the
$\langle \overline R\,R\rangle$-propagators:
$$
{\psset{unit=1}
\begin{pspicture}(-7,-1.2)(7,1.2)
\newcommand{\PATTERNONE}{%
{\psset{unit=1}
\psarc[linecolor=white, linestyle=solid, linewidth=10pt](0,0){1}{30}{150}
\psarc[linecolor=green, linestyle=dotted, linewidth=2pt](0,0){.8}{30}{150}
\psarc[linecolor=blue, linestyle=solid, linewidth=1pt](0,0){1.2}{30}{150}
\rput(-1.1,1){\makebox(0,0)[cc]{\hbox{{$\Lambda$}}}}
\rput(1.1,1.05){\makebox(0,0)[cc]{\hbox{{$\overline\Lambda$}}}}
\rput(0,.5){\makebox(0,0)[cc]{\hbox{{$\beta$}}}}
}
}
\newcommand{\PATTERNTWO}{%
{\psset{unit=1}
\psarc[linecolor=white, linestyle=solid, linewidth=10pt](0,0){1}{210}{330}
\psarc[linecolor=green, linestyle=dotted, linewidth=2pt](0,0){.8}{210}{330}
\psarc[linecolor=red, linestyle=dashed, linewidth=1.5pt](0,0){1.2}{210}{330}
}
}
\newcommand{\PATTERNTHREE}{%
{\psset{unit=1}
\pcline[linecolor=white, linestyle=solid, linewidth=10pt](0,0)(1,0)
\pcline[linecolor=blue, linestyle=solid, linewidth=1pt](0,0.2)(1,0.2)
\pcline[linecolor=red, linestyle=dashed, linewidth=1.5pt](0,-.2)(1,-.2)
\rput(1,0){\pscircle[linecolor=black, fillstyle=solid, fillcolor=white]{.2}}
\rput(0,0){\pscircle*{.2}}
}
}
\rput(-2.86,0.5){\PATTERNTWO}
\rput(1.86,0.5){\PATTERNTWO}
\rput(4.59,0.5){\PATTERNTWO}
\rput(-2.86,-0.5){\PATTERNONE}
\rput(1.86,-0.5){\PATTERNONE}
\rput(4.59,-0.5){\PATTERNONE}
\rput(-6.73,0){\PATTERNTHREE}
\rput(-4.73,0){\PATTERNTHREE}
\rput(-2,0){\PATTERNTHREE}
\rput(0,0){\PATTERNTHREE}
\rput(2.73,0){\PATTERNTHREE}
\rput(5.46,0){\PATTERNTHREE}
\rput(-5.23,0){\makebox(0,0)[cc]{\hbox{{$+$}}}}
\rput(-.5,0){\makebox(0,0)[cc]{\hbox{{$+$}}}}
\rput(7,0){\makebox(0,0)[lc]{\hbox{{$+\cdots$}}}}
\rput(-6.73,-0.3){\makebox(0,0)[ct]{\hbox{{$\overline R$}}}}
\rput(-4.73,-0.3){\makebox(0,0)[ct]{\hbox{{$\overline R$}}}}
\rput(0,-0.3){\makebox(0,0)[ct]{\hbox{{$\overline R$}}}}
\rput(-5.73,-0.35){\makebox(0,0)[ct]{\hbox{{$R$}}}}
\rput(-1,-0.35){\makebox(0,0)[ct]{\hbox{{$R$}}}}
\rput(6.46,-0.35){\makebox(0,0)[ct]{\hbox{{$R$}}}}
\end{pspicture}
}
$$
so that the propagator becomes
$$
\langle \overline R\,R\rangle\sim \frac{1}{N}\frac{\delta_{i_1,i_2}\delta_{k_1,k_2}}{1-\beta|\lambda_{i_1}|^2}
$$
and the corresponding quadratic form gets an external field addition:
\be
-N\tr[\overline R R(1-\beta |\Lambda|^2)].
\ee
The new interaction vertex arises from the dotted cycles of length four:
$$
{\psset{unit=1}
\begin{pspicture}(-2,-2)(4,2)
\newcommand{\PATTERNTHREE}{%
{\psset{unit=1}
\pcline[linecolor=white, linestyle=solid, linewidth=10pt](0,0)(0.8,0)
\pcline[linecolor=blue, linestyle=solid, linewidth=1pt](0,0.2)(0.8,0.2)
\pcline[linecolor=red, linestyle=dashed, linewidth=1.5pt](0,-.2)(0.8,-.2)
}
}
\rput{45}(0,0){
\rput(1,0){\PATTERNTHREE
}
}
\rput{225}(0,0){
\rput(1,0){
\PATTERNTHREE
}
}
\rput{315}(0,0){
\rput(-1.8,0){\PATTERNTHREE
}
}
\rput{135}(0,0){
\rput(-1.8,0){\PATTERNTHREE
}
}
\pscircle[linecolor=white, fillstyle=solid, fillcolor=white](0,0){1.2}
\pscircle[linecolor=green, linestyle=dotted, linewidth=2pt](0,0){.8}
\psarc[linecolor=blue, linestyle=solid, linewidth=1pt](0,0){1.2}{55}{125}
\psarc[linecolor=blue, linestyle=solid, linewidth=1pt](0,0){1.2}{235}{305}
\psarc[linecolor=red, linestyle=dashed, linewidth=1.5pt](0,0){1.2}{145}{215}
\psarc[linecolor=red, linestyle=dashed, linewidth=1.5pt](0,0){1.2}{-35}{35}
\rput(-0.5,1.2){\makebox(0,0)[rb]{\hbox{{$\Lambda$}}}}
\rput(0.5,1.2){\makebox(0,0)[lb]{\hbox{{$\overline\Lambda$}}}}
\rput(0.5,-1.2){\makebox(0,0)[lt]{\hbox{{$\Lambda$}}}}
\rput(-0.5,-1.2){\makebox(0,0)[rt]{\hbox{{$\overline\Lambda$}}}}
\rput(0,0){\makebox(0,0)[cc]{\hbox{{$\beta$}}}}
\rput{45}(0,0){
\rput(1,0){\pscircle*{.2}}
}
\rput{225}(0,0){
\rput(1,0){\pscircle*{.2}}
}
\rput{135}(0,0){
\rput(1,0){\pscircle[linecolor=black, fillstyle=solid, fillcolor=white]{.2}}
}
\rput{315}(0,0){
\rput(1,0){\pscircle[linecolor=black, fillstyle=solid, fillcolor=white]{.2}}
}
\rput(1.6,0){\makebox(0,0)[lc]{\hbox{{$\sim \frac 12 N\beta\tr[\overline R R \Lambda\overline \Lambda\,\overline R R \Lambda\overline \Lambda]$}}}}
\end{pspicture}
}
$$
where the factor $1/2$ takes into account the symmetry of the four-cycle.

We therefore have that the generating function ${\mathcal F}$ is the logarithm of the integral
\be
\int DR\,D\overline R\, e^{N\tr[-\overline R R(1-\beta |\Lambda|^2)+\frac12 \beta
\overline R R|\Lambda|^2\overline R R|\Lambda|^2]},
\label{int-clean1}
\ee
where we integrate over rectangular complex $(\gamma N\times \alpha N)$-matrices $R$.
We first rescale the integration variable $R\to R\Lambda$, which results in the integral
\be
\prod_{i=1}^{\alpha N}|\lambda_i|^{-2\gamma N}
\int DR\,D\overline R\, e^{N\tr[-\overline R R( |\Lambda|^{-2}-\beta ) + \frac12 \beta
\overline R R\overline R R]}.
\label{int-clean2}
\ee
Performing now the same chain of transformations as in Sec.~\ref{s:model}, we obtain eventually that
integral (\ref{int-clean2}) is equivalent to the Hermitian one-matrix model integral
\be
\prod_{i=1}^{\alpha N}|\lambda_i|^{-2\gamma N}
\int_{\alpha N\times \alpha N} DH_{\ge 0}
e^{N\tr[-H(\Lambda^{-2}-\beta)+(\gamma-\alpha)\log H +\frac 12 \beta H^2]}.
\label{int-clean3}
\ee

\begin{lm}\label{lm:gen-Groth-clean}
The generating function for clean Belyi fat graphs ((\ref{gen-fun}) with ramification profiles $(2,\dots,2,1,\dots,1)$ at the point $1$) is the matrix-model integral (\ref{int-clean3}). This matrix-model integral is the (original) Kontsevich--Penner matrix model \cite{ChM}, \cite{ChM2}.
\end{lm}

\begin{remark}
If we demand the ramification profile at the point $1$ to be just $(2,2,\dots,2)$ (no dotted two-cycles are allowed),
then in order to obtain the corresponding generating function we must merely replace $\Lambda^{-2}-\beta$ by $\Lambda^{-2}$ in (\ref{int-clean3}).
\end{remark}

From now on, for simplicity, we restrict ourselves to the case of ramification profile $(2,2,\dots,2)$ at the
point $1$.

\subsection{Solving integral (\ref{int-clean3})}
That the Kontsevich--Penner matrix model integral (\ref{int-clean3}) is equivalent to the Hermitian one-matrix model
integral (\ref{1MM}) is well known. This equivalence was established using the Virasoro constraints in \cite{ChM2}
or using explicit determinant relations in \cite{MMM}. We recall here the logic of \cite{MMM}.

We begin with the standard eigenvalue representation for integral (\ref{1MM}),
\be
\int {dy_1\dots dy_M}\,[\Delta(y)]^2e^{-\sum_{k=1}^\infty \sum_{i=1}^M \xi_k y_i^k}
\label{int-1MM}
\ee
in which we again perform the Miwa change of variables with the Gaussian shift,
\be
\xi_k=\frac 1k \sum_{j=1}^N \frac 1{\mu_j^k}+\frac 12 \delta_{k,2}.
\label{Miwa}
\ee
Summing up the terms in the exponential into logarithms, we transform integral (\ref{int-1MM}) to the form
$$
\int {dy_1\dots dy_M}\,[\Delta(y)]^2\prod_{i=1}^M\prod_{j=1}^N(\mu_j-y_i)\prod_{j=1}^N \mu_j^{-M}
e^{-\frac12\sum_{i=1}^M h_i^2}.
$$
We now use that $\Delta(y)\prod_{i=1}^M\prod_{j=1}^N(\mu_j-y_i)=\Delta(y,\mu)/\Delta(\mu)$, where
$\Delta(y,\mu)$ is the Vandermonde determinant of the set of variables $y_i$ and $\mu_j$, write each of the
determinants $\Delta(y,\mu)$ and $\Delta(y)$ as determinants of the Hermitian polynomials $H_s(x)$, where
$s$ ranges from $0$ to $M+N-1$ and $x$ are either $y_i$ or $\mu_j$ in the first determinant and $s$ ranges
from $0$ to $M-1$ and $x$ are $y_i$ in the second determinant. Because the Hermitian polynomials are orthogonal
with the measure $e^{-\frac 12 x^2}$, we can integrate out all the $y$-variables; the remaining expression will be
the determinant of the $(N\times N)$-matrix $\| H_{M+{j_1}-1}(\mu_{j_2})\|$, $j_1,j_2=1,\dots,N$, and the
original integral (\ref{int-1MM}) thus takes the form
\be
\prod_{j=1}^N \mu_j^{-M} \frac1{\Delta(\mu)}\left|\begin{array}{cccc}
H_M(\mu_1) & H_M(\mu_2)& \dots& H_M(\mu_N)\cr
H_{M+1}(\mu_1) & H_{M+1}(\mu_2)& \dots& H_{M+1}(\mu_N)\cr
\vdots &\vdots &\dots &\vdots\cr
H_{M+N-1}(\mu_1) & H_{M+N-1}(\mu_2)& \dots& H_{M+N-1}(\mu_N)
\end{array}
\right|
\ee
On the other hand, we obtain the same ratio of determinants multiplied by $e^{-\frac 12 \sum_j \mu_j^2}$
if we consider the $N$-fold integral
\be
\int dx_1\dots dx_N\,\frac{\Delta(x)}{\Delta(\mu)}\prod_{j=1}^N x_j^M e^{\sum_{j=1}^N (x_j\mu_j+\frac 12 x_j^2)}
\label{extern}
\ee
because $\int dx\, x^s e^{x\mu+\frac 12x^2}=e^{-\frac 12\mu^2}H_s(\mu)$. Expression (\ref{extern}) is nothing but the
Kontsevich--Penner integral, so we obtain the relation between two matrix integrals of {\it different} sizes:
\be
\int\limits_{N\times N} DX\,e^{\tr[X\mu +\frac 12 X^2 + M\log X]}=\prod_{j=1}^M\bigl[\mu_j^M
e^{-\frac12 \mu_j^2}\bigr]\int\limits_{M\times M}DY\,e^{-\sum\limits_{k=1}^\infty\xi_k\tr Y^k},
\ \xi_k={\textstyle\frac1k \sum\limits_{j=1}^N\frac1{\mu_j^k}+\frac 12\delta_{k,2}}.
\label{5.9}
\ee

After a simple algebra, we come to the following lemma.
\begin{lm}\label{lm:clean}
The generating function (\ref{gen-fun})
for the {\em clean} Belyi morphisms with the ramification profile $(2,2,\dots,2)$ at the
point $1$ is given by the following Hermitian one-matrix model integral for $\gamma-\alpha\simeq O(1)$:
\be
{\mathcal Z}[t;\gamma,\beta]=
\prod_{i=1}^{\alpha N}|\lambda_i|^{-2\gamma N}
\int\limits_{M\times M}DY\,e^{-\sum\limits_{k=1}^{\infty}\frac{t_k}{k}(-1)^k\tr Y^k-\frac{N}{2\beta}\tr Y^2},
\ t_k={\textstyle\sum\limits_{i=1}^{\alpha N}}\lambda_i^{2k},\ M=(\gamma-\alpha)N.
\label{clean1}
\ee
Because this integral is also equivalent to Kontsevich--Penner matrix model (\ref{int-clean3}) (with the external
field term $\Lambda^{-2}$ instead of $\Lambda^{-2}-\beta$), it also belongs to the GKM class thus being a tau
function of the KP hierarchy.
\end{lm}

\begin{remark}\label{rm:no-limit-no}
Note again that the above correspondence is valid only in the $1/N$ asymptotic expansion and only when
$\gamma-\alpha\simeq O(1)$. If $\gamma-\alpha\lesssim O(1/N)$ the above correspondence fails because in this case
we must take into account that we integrate in formula (\ref{int-clean3}) over positive definite
matrices, contrary to formula (\ref{5.9}) in which no restriction on integration domain is assumed. So, again,
the case $\gamma=\alpha$ is special and must be treated separately.
\end{remark}

\section{A general case of two-profile Belyi morphisms}\label{s:general}
\setcounter{equation}{0}

Combining the techniques of Secs.~\ref{s:model} and~\ref{s:CEO} we now address the most general case of
Belyi morphisms with the given profiles at {\it two} branching points: infinity and $1$. We take these profiles into
account in two different ways: at infinity we, as in Sec.~\ref{s:CEO}, introduce the times $t_m$ responsible for
the profile whereas the times at $1$ will be taken into account by introducing, as in Sec.~\ref{s:model},
the external field $\Lambda$ with
\be
{\mathfrak t}_s=\tr\bigl[(\Lambda\overline \Lambda)^s\bigr]=\sum_{k=1}^{\gamma N}|\lambda_k|^{2s}.
\label{tau-profile}
\ee
We then have the following statement
\begin{lm}
The generating function
\be
{\mathcal F}[\{t_1,t_2,\dots\},\{{\mathfrak t}_1,{\mathfrak t}_2,\dots\},\beta;N]=\sum_{\Gamma}\frac{1}{|\hbox{Aut\,}\Gamma|}
N^{2-2g}\beta^{n_2}\prod_{i=1}^{n_1}t_{r_i}\prod_{k=1}^{n_3}{\mathfrak t}_{s_k}
\label{gen-fun-gen}
\ee
of Belyi morphisms in which we have two sets of ramification profiles: $\{t_{r_1},\dots,t_{r_{n_1}}\}$ at infinity
and $\{{\mathfrak t}_{s_1},\dots,{\mathfrak t}_{s_{n_3}}\}$ at $1$ is given by the integral over complex rectangular
$(\beta N\times \gamma N)$-matrices $B,\overline B$:
\be
{\mathcal Z}[t,{\mathfrak t}]:=e^{{\mathcal F}[\{t\},\{{\mathfrak t}\},\beta;N]}
=\int_{\gamma N\times \beta N}DB\,D\overline B\,
e^{-N\tr[B\overline B]+N\sum_{m=1}^\infty \frac 1m t_m\tr [(B\overline B\,\overline\Lambda\Lambda)^m]},
\label{tt1}
\ee
where the times ${\mathfrak t}_s$ are given by (\ref{tau-profile}).
\end{lm}

Performing the same operation as in (\ref{new1})--(\ref{new2}), we obtain that integral (\ref{tt1})
is equal to the integral over Hermitian positive definite $(\gamma N\times \gamma N)$-matrix $X$
with the external matrix field $\tilde\Lambda=|\Lambda|^{-2}$:
\be
{\mathcal Z}[t,\tau]=\prod_{k=1}^{\gamma N}|\lambda_k|^{-2\beta N}
\int\limits_{\gamma N\times \gamma N} DX_{\ge 0}
e^{N\tr\Bigl[-X|\Lambda|^{-2}+ \sum\limits_{m=1}^\infty \frac {t_m}m X^m+(\beta-\gamma)\log X\Bigr]},
\label{tt2}
\ee
Integral (\ref{tt2}) is again a GKM integral \cite{MMM}; after integration
over eigenvalues $x_k$ of the matrix $X$ it takes the form of the ratio of two determinants,
\be
{\mathcal Z}[t,\tau]=\prod_{k=1}^{\gamma N}|\lambda_k|^{-2\beta N}
\frac{\Bigl\| \frac{\partial^{k_1-1}}{\partial {\tilde\lambda}_{k_2}^{k_1-1}}f(\tilde\lambda_{k_2})
\Bigr\|_{k_1,k_2=1}^{\gamma N}}{\Delta(\tilde\lambda)},
\label{tt3}
\ee
where
\be
f(\tilde\lambda)=\int_{0}^{\infty}x^{N(\beta-\gamma)}
e^{-Nx\tilde\lambda+N\sum\limits_{m=1}^\infty \frac {t_m}m x^m}.
\label{tt4}
\ee
Because any GKM integral (in the proper normalization) is a $\tau$-function of the
KP hierarchy, and for a model with the logarithmic term in the potential it was demonstrated in
\cite{MMS}, we immediately obtain the following theorem.

\begin{theorem}
The exponential $e^{{\mathcal F}[\{t\},\{{\mathfrak t}\},\gamma;N]}$
of generating function (\ref{gen-fun-gen}) modulo the
normalization factor $\prod_{k=1}^{\gamma N}|\lambda_k|^{-2\beta N}$
is a $\tau$-function of the KP hierarchy (that is, it satisfies the bilinear Hirota relations)
in times ${\mathfrak t}_s$ given by (\ref{tau-profile}).
\end{theorem}

\section{Conclusion}
We have proved that generating functions for numbers of three different types of Belyi morphisms are free
energies of special matrix models all of which are in the GKM class thus being tau functions of the KP hierarchy.
Besides this, it is interesting to establish other relations between, say, generating function (\ref{gen-fun})
for clean Belyi morphisms and the free energy of the Kontsevich--Penner matrix model, which is known (see
\cite{Ch95},\cite{Norbury},\cite{DoNor}) to be related to the numbers of integer points in moduli spaces
${\mathcal M}_{g,n}$ of curves of genus $g$ with $n$ holes with fixed (integer) perimeters; the very same model
is also related \cite{Ch95} by a canonical transformation to two copies of the Kontsevich matrix model expressed
in times related to the discretization of the moduli spaces ${\mathcal M}_{g,n}$. It is tempting to find possible
relations between these discretizations, cut-and-join operators of \cite{Zograf},
and Hodge integrals of \cite{Kazarian}.

Of course, the possibility of using GKM techniques when studying enumeration problems for Belyi morphisms deserves
more detailed studies; we consider this note a first step in exploring this perspective field of knowledge.

It is also interesting to clarify the role of cut-and-join operators of \cite{Kazarian} and \cite{Zograf}
in the matrix-model context. After this text was completed, an interesting paper \cite{AMMN} extending the
formalism of cut-and-join operators to the case of generalized Hurwitz numbers has appeared.

\section*{Acknowledgments}
The authors thank Maxim Kazarian, Andrei Mironov, and Petr Zograf for the useful discussion.

The authors acknowledge support from the ERC Advance Grant 291092 ``Exploring the Quantum Universe'' (EQU).
J.A. acknowledges support of the FNU, the Free Danish Research Council, from the grant ``Quantum gravity and the
role of black holes.''  The work of J.A. was supported in part by Perimeter Institute of Theoretical Physics. Research of
Perimeter Institute is supported by the Government of Canada through Industry Canada and by Province of Ontario through
the Ministry of Economic Development and Innovation. 
The work of L.Ch. was supported by the Russian Foundation for Basic Research
(Grant Nos. 14-01-00860-a and 13-01-12405-ofi-m) and by the Program Mathematical Methods for Nonlinear Dynamics.

\def\thetheorem{\Alph{section}.\arabic{theorem}}
\def\theprop{\Alph{section}.\arabic{prop}}
\def\thelemma{\Alph{section}.\arabic{lm}}
\def\thecor{\Alph{section}.\arabic{cor}}
\def\theexam{\Alph{section}.\arabic{exam}}
\def\theremark{\Alph{section}.\arabic{remark}}
\def\theequation{\Alph{section}.\arabic{equation}}

\setcounter{section}{0}

\appendix{Deriving the Jacobian of transformation (\ref{RR2UVM})}\label{se:notation}
\setcounter{equation}{0}

The invariant measure $DU\,DV$ in the vicinity of the unity becomes $DH\,D\tilde H\, DP\,D\overline P$. For
$d\overline R_{i,k}$ we then obtain
\be
d\overline R_{i,k}=\Bigl\{d\overline m_{i}\delta_{i,k}+ idH_{i,k}\overline m_k +im_i d\tilde H_{i,k},\ k\le\alpha N
\Bigm| \overline m_i dP_{i,k-\alpha N}, \ k>\alpha N\Bigr\}.
\label{A0}
\ee
The elements $dm_{i}$ appear only for $i=k$ with the unit factor, so we have to calculate only ``non-diagonal''
differentials $DR\,D\overline R$. For $i<k\le\alpha N$ we have:
\be
\begin{array}{ll}
d\overline R_{i,k}=idH_{i,k}\overline m_k +i\overline m_i d\tilde H_{i,k},&
d\overline R_{k,i}=idH^*_{i,k}\overline m_i +i \overline m_k d\tilde H^*_{i,k},\phantom{\Bigl|}\cr
dR_{k,i}=-idH_{i,k} m_i-im_k d\tilde H_{i,k}, &
dR_{i,k}=-idH^*_{i,k} m_k-i m_id\tilde H^*_{i,k}.\phantom{\Bigl|}
\end{array}
\ee
Combining the columns in these relations, we obtain
\be
\begin{array}{l}
d\overline R_{i,k}\wedge dR_{k,i}
=dH_{i,k}\wedge d\tilde H_{i,k}[m_k\overline m_k-m_i\overline m_i],\phantom{\Bigl|}\cr
d\overline R_{k,i}\wedge dR_{i,k}=dH^*_{i,k}\wedge d\tilde H^*_{i,k}[m_i\overline m_i-m_k\overline m_k],
\end{array}
\ 1\le i<k\le \alpha N,
\label{A1}
\ee
and we obtain that
\be
\mathop{\hbox{\large$\wedge$}}\limits_{i,k=1}^{\alpha N}d\overline R_{i,k}\wedge dR_{k,i}=DH\wedge D\tilde H\wedge
\prod_{i=1}^{\alpha N}dm_i\wedge d\overline m_i \prod_{1\le i<k\le \alpha N}\bigl[|m_i|^2-|m_k|^2\bigr]^2.
\ee
For the remaining part we merely obtain from (\ref{A0}) that
\be
\mathop{\hbox{\large$\wedge$}}\limits_{{i=1,\dots,\alpha N\atop k=\alpha N+1,\dots,\gamma N}}
d\overline R_{i,k}\wedge dR_{k,i}
=DP\wedge D\overline P \prod_{i=1}^{\alpha N}|m_i|^{2(\gamma-\alpha)N},
\ee
so we finally obtain formula (\ref{Jac}) for the Jacobian of transformation (\ref{RR2UVM}).

\end{document}